\documentclass[10pt]{article}

\usepackage{amsmath,amsfonts,amsthm,amssymb,graphicx,mathdots}
\usepackage[all,2cell,ps]{xy}

\theoremstyle{plain}
\newtheorem{thm}{Theorem}[section]
\newtheorem{lem}[thm]{Lemma}
\newtheorem{prop}[thm]{Proposition}
\newtheorem{cor}[thm]{Corollary}

\theoremstyle{definition}

\newtheorem*{qtn}{Question}
\newtheorem*{rem}{Remark}

\theoremstyle{remark}

\newcommand{\ep}{\epsilon}

\DeclareMathOperator{\U}{U}

\newcommand{\al}{\alpha}
\newcommand{\F}{\mathbb F}
\newcommand{\Z}{\mathbb Z}

\newcommand{\R}{\mathbb R}

\newcommand{\Q}{\mathbb Q}

\newcommand{\gam}{\gamma}

\newcommand{\lam}{\lambda}

\newcommand{\A}{\mathbb{A}}
\newcommand{\C}{\mathbb C}
\newcommand{\Hy}{\mathbf H}

\newcommand{\Gam}{\Gamma}

\newcommand{\conj}{\overline}

\newcommand{\mc}{\mathcal}

\DeclareMathOperator{\SL}{SL}
\DeclareMathOperator{\PSL}{PSL}

\DeclareMathOperator{\PGL}{PGL}

\DeclareMathOperator{\SO}{SO}

\DeclareMathOperator{\PU}{PU}
\DeclareMathOperator{\SU}{SU}

\DeclareMathOperator{\real}{Re}

\newcommand{\mf}{\mathfrak}
\DeclareMathOperator{\Hom}{Hom}

\DeclareMathOperator{\Aut}{Aut}

\newcommand{\ssm}{\smallsetminus}

\DeclareMathOperator{\Id}{Id}

\newenvironment{pf}{\begin{proof}}{\end{proof}}

\newenvironment{mat}{\left(\begin{matrix}}{\end{matrix}\right)}

\title{Cusps of Picard modular surfaces}

\author{Matthew Stover \\ \small{University of Texas at Austin}\\ \small{\textsf{mstover@math.utexas.edu}}}

\date{\today}

\begin{document}

\maketitle

\begin{abstract}
We determine the number of cusps of minimal Picard modular surfaces. The proof also counts cusps of other Picard modular surfaces of arithmetic interest. Consequently, for each $N > 0$ there are finitely many commensurability classes of nonuniform arithmetic lattices in $\SU(2,1)$ that contain an $N$-cusped surface. We also discuss a higher-rank analogue.
\end{abstract}



\section{Introduction}\label{intro}


To study the cusps of locally symmetric spaces, it is most natural to begin with the minimal elements of each commensurability class. From there, one can, ostensibly, use the combinatorics of covering spaces to study the rest of the commensurability class. The purpose of this paper is to study the cusps of minimal arithmetic complex hyperbolic 2-orbifolds. Equivalently, we study cusps of maximal nonuniform arithmetic lattices in $\PU(2, 1)$. These algebraic surfaces are often called Picard modular surfaces, and their fundamental groups Picard modular groups.

We begin with two remarks. First, by Mostow--Prasad rigidity, when $\Gam$ is an irreducible lattice in any semisimple group except $\SL_2(\R)$, the corresponding locally symmetric quotient $M$ is uniquely determined by $\Gam$. Therefore, there is no ambiguity in saying `the cusps of $\Gam$' as opposed to `the cusps of $M$'. For $\SL_2(\R)$, there is an easy counterexample: the free group on two generators is the fundamental group of both the thrice-punctured sphere and once-punctured torus. Second, it is a consequence of work of Borel and the arithmeticity theorem of Margulis that a lattice $\Gam$ in a semisimple Lie group is arithmetic if and only if there are infinitely many distinct isomorphism classes of maximal lattices in its commensurability class. In fact, $\Gam$ is nonarithmetic if and only if there is a \emph{unique} maximal element in its commensurability class.

If $\Gam < \PU(2,1)$ is a maximal arithmetic lattice, it is best described using Bruhat--Tits theory. See $\S$\ref{arithmetic}. The preceding section, $\S$\ref{prelims}, is a brief review of complex hyperbolic geometry and some of the algebraic tools we use. The goal of $\S$\ref{cusps} is to count cusps of minimal Picard modular surfaces. The results are difficult to state without the notation of $\S$\ref{arithmetic}, so we refer to $\S$\ref{cusps} for precise statements. The main technical result is Theorem \ref{maximal count}. The formula for a minimal surface is roughly the following.


\begin{thm}\label{intro cusp count}
Let $M$ be a minimal Picard modular surface, defined via a hermitian form over the imaginary quadratic field $k$. Then there exists a nonnegative integer $n$, depending explicitly on $M$, such that $M$ has
\[
3^n \frac{h_k}{h_{k, 3}}
\]
cusps, where $h_k$ is the class number of $k$ and $h_{k, 3}$ is the order of the three-primary part of the class group of $k$.
\end{thm}


In short, the proof involves counting the number of cusps for particular congruence subgroups of a given lattice, then using these counts to determine the number of cusps for the minimal surface. These congruence subgroups are natural generalizations of Hecke's subgroups $\Gam_0(N)$ in $\PSL_2(\Z)$ (cf.\ the groups considered in \cite{Chinburg--Long--Reid}). Using Bruhat--Tits theory, one is then able to move around the commensurability class and keep track of the action on cusps. The proof is inspired equally by Zink's study of cusps of lattices in unitary groups \cite{Zink} and Chinburg--Long--Reid's paper on cusps of minimal arithmetic hyperbolic 3-orbifolds \cite{Chinburg--Long--Reid}.

Our primary application of Theorem \ref{intro cusp count} is to the classification, up to commensurability, of noncompact finite volume complex hyperbolic orbifolds. The following, which we prove at the end of $\S$\ref{cusps}, utilizes work of Ellenberg--Venkatesh \cite{Ellenberg--Venkatesh} on torsion in the ideal class group.


\begin{thm}\label{intro cusp commensurability}
For any $N > 0$, there are only finitely many commensurability classes of Picard modular surfaces containing an element with $N$ cusps.
\end{thm}


It remains unknown whether or not there are infinitely many commensurability classes of nonarithmetic complex hyperbolic 2-orbifolds. For hyperbolic 3-space, Chinburg--Long--Reid \cite{Chinburg--Long--Reid} proved the analogue to Theorem \ref{intro cusp commensurability}, but there are infinitely many commensurability classes of non-commensurable hyperbolic knot complements (e.g., twist knots), so arithmetic assumption is crucial.

We also prove a higher-rank version of Theorems \ref{intro cusp count} and \ref{intro cusp commensurability}. Let $\Gam < \SU(r + 1, r)$ be a maximal nonuniform lattice with $r > 1$. By the monumental work of Margulis, $\Gam$ is automatically arithmetic. All nonuniform arithmetic lattices in $\SU(2, 1)$ arise from hermitian forms in $k^3$, where $k$ is any imaginary quadratic field. In higher rank, there are other methods of producing nonuniform arithmetic lattices, all of which arise from hermitian forms over division algebras with involution of second kind. We say that a lattice is of \emph{simple type} if it arises from a hermitian form over an imaginary quadratic field. As in \cite{Zink}, the methods of Theorem \ref{intro cusp count} extend to higher rank lattices of simple type.


\begin{thm}\label{intro higher rank cusp count}
Let $\Gam < \SU(r + 1, r)$ be a maximal nonuniform lattice of simple type, and set $q = 2 r + 1$. Then there exists a number $n$, depending explicitly on $\Gam$, so that $\Gam$ has
\[
q^n \frac{h_k^r}{h_{k, q}}
\]
cusps, where $h_k$ is the class number of $k$ and $h_{k, q}$ is the order of the subgroup of the class group of $k$ consisting of elements with order dividing $q$.
\end{thm}


Since $\SU(2, 1)$ has $\R$-rank one, cusps correspond to topological ends of the corresponding orbifold. In higher rank, `cusps' correspond to maximal chains $P_0 \subset P_1 \subset \cdots \subset P_r$ of parabolic subgroups. For geometric interpretations, see \cite{Borel--Ji}. When $q$ is prime, all nonuniform lattices are of simple type. This allows us to give a full classification of the commensurability classes of locally symmetric spaces with $\SU(r + 1, r)$-geometry containing a one-cusped element.


\begin{thm}\label{intro higher rank commensurability}
For $r > 1$ and $2 r + 1$ prime, there are only finitely many commensurability classes of locally symmetric spaces with $\SU(r + 1, r)$-geometry containing a one-cusped element. Furthermore, assuming the Generalized Riemann Hypothesis, there are only finitely many commensurability classes containing an $N$-cusped element for any $N > 0$.
\end{thm}


We also show how, for each $r \geq 1$, one can construct infinitely many homeomorphism classes of one-cusped orbifolds with $\SU(r + 1, r)$-geometry. However, none of the spaces we build are manifolds. The following seems to be a significant gap in our understanding of complex hyperbolic manifolds.


\begin{qtn}
Does there exist a one-cusped complex hyperbolic 2-manifold?
\end{qtn}


Restrictions on the possible cusp cross-sections are given by unpublished work of Neumann--Reid and Kamishima \cite{Kamishima}. There are one-cusped hyperbolic manifolds commensurable with the modular curve, and, though the figure-eight knot is the unique arithmetic knot complement \cite{Reid}, several Bianchi groups $\PSL_2(\mc{O}_d)$ contain a one-cusped manifold. See also \cite{Petersen}. The corresponding question for higher dimensional hyperbolic manifolds also appears to be unsolved.


\section{Preliminaries}\label{prelims}



\subsection{Complex hyperbolic geometry}\label{geometry}


Let $V$ be a three-dimensional complex vector space equipped with a nondegenerate hermitian form $h$ of signature $(2,1)$. In this paper, $h$ will typically be the form with matrix
\[
h_0 = \begin{mat}
0 &
0 &
1 \\

0 &
-1 &
0 \\

1 &
0 &
0
\end{mat}.
\]
The \emph{complex hyperbolic plane} is $\mathbf{H}_\C^2 = \mathbb{P}(V_-) \subset \mathbb{P}^2$, where $V_-$ is the set of $h$-positive vectors in $V$. With the metric associated to $h$, the biholomorphic isometry group of $\mathbf{H}_\C^2$ is $\PU(2,1)$.

The \emph{ideal boundary} $\partial_\infty \mathbf{H}_\C^2$ of complex hyperbolic space is the image in $\mathbb{P}^2$ of the $h$-isotropic vectors $V_0$. For $h_0$, $\mathbf{H}_\C^2$ has homogeneous coordinates
\[
\left\{ \left[ \begin{matrix}
z_1 \\

z_2 \\

z_3
\end{matrix} \right]\ :\ 2 \real(z_1\conj{z}_3) - |z_2|^2 > 0 \right\}
\]
and $\partial_\infty \mathbf{H}_\C^2$ is
\[
\left\{ \left[ \begin{matrix}
z_1 \\

z_2 \\

z_3
\end{matrix} \right]\ :\ 2 \real(z_1\conj{z}_3) - |z_2|^2 = 0 \right\}.
\]

A \emph{complex hyperbolic 2-orbifold} is $\mathbf{H}_\C^2 / \Gam$, where $\Gam < \PU(2,1)$ is a discrete subgroup of isometries such that $\mathbf{H}_\C^2 / \Gam$ has finite volume. Isometries fall into three categories based upon their action on $\mathbf{H}_\C^2$ and its ideal boundary. An isometry is \emph{elliptic} if it has a fixed point in the complex hyperbolic plane, \emph{loxodromic} if it is not elliptic and has two fixed points on the ideal boundary, and is \emph{parabolic} if it is not elliptic and fixes exactly one point on the ideal boundary.

A lattice $\Gam$ is called \emph{uniform} if $\mathbf{H}_\C^2 / \Gam$, or equivalently $\PU(2,1) / \Gam$, is compact. Otherwise it is \emph{nonuniform} and $\mathbf{H}_\C^2 / \Gam$ has a finite number of topological ends. A lattice is nonuniform if and only if it contains parabolic elements.

The \emph{cusp set} of $\Gam$ is the set of all fixed points of parabolic elements in $\Gam$. There is a natural $\Gam$-action on the cusp set, and a \emph{cusp} of $\Gam$ is a $\Gam$-equivalence class of fixed points of parabolic elements in $\Gam$. The topological ends of $\mathbf{H}_\C^2 / \Gam$ are in one-to-one correspondence with the cusps of $\Gam$. As mentioned in the introduction, by Mostow--Prasad rigidity, the number of cusps of $\mathbf{H}_\C^2 / \Gam$ is an invariant of the lattice. Therefore, we will alternate freely between saying `the cusps of $\Gam$' and `the cusps of $\mathbf{H}_\C^2 / \Gam$' throughout.


\subsection{Algebra of hermitian forms}\label{algebra}


Let $V$ be a vector space of dimension $n$ over the number field $k$. If $\mc{O}_k$ is the ring of integers of $k$, a \emph{lattice} $\mc{L} \subset V$ is an $\mc{O}_k$-submodule such that $\mc{L} \otimes_{\mc{O}_k} k = V$. Since $\mc{O}_k$ is a Noetherian ring, there is an isomorphism of $\mc{O}_k$-modules
\[
\mc{L} \cong_{\mc{O}_k} I_1 \oplus \cdots \oplus I_n,
\]
where each $I_j$ is a fractional ideal of $k$. This decomposition is determined up to isomorphism by the \emph{ideal class} $\mathrm{cl}(\mc{L})$, which is the class of the product $I_1 \cdots I_n$ in the ideal class group of $k$. That is, two lattices $\mc{L}_1$ and $\mc{L}_2$ are isomorphic as $\mc{O}_k$-modules if and only if they have the same dimension and $\mathrm{cl}(\mc{L}_1) = \mathrm{cl}(\mc{L}_2)$.

Now, suppose that $k$ is an imaginary quadratic field and that $h$ is a nondegenerate hermitian form on $V$. Then $h$ is determined up to isomorphism by $\det(h)$, considered as an element of $\Q^\times / \mathrm{N}_{k / \Q}(k^\times)$ \cite[Chapter 10]{Scharlau}. A \emph{hermitian lattice} is a pair $(\mc{L}, h)$, where $\mc{L}$ is a lattice in $V$ and $h$ a hermitian form. Hermitian lattices are then classified up to isomorphism by $\mathrm{cl}(\mc{L})$ and $\det(h)$.

We say that a hermitian lattice $(\mc{L}, h)$ is \emph{unimodular} if
\[
\mc{L} = \{ x \in V\ :\ h(x, \mc{L}) \subset \mc{O}_k \}.
\]
The lattice $h_0$ from $\S$\ref{geometry} is unimodular with respect to the standard lattice $\mc{O}_k^3$ for that basis. More generally, if $I$ is a fractional ideal of $k$, we call $(\mc{L}, h)$ $I$-\emph{modular} if
\[
\mc{L} = \{ x \in V\ :\ h(x, \mc{L}) \subset I \}.
\]
If $V$ is odd-dimensional and $(\mc{L}, h)$ is unimodular, then $\mc{L}$ contains an odd vector, so given a reduced element $x \in \mc{L}$, there exists $y \in \mc{L}$ such that $h(x, y) = 1$. This fact will be crucial in what follows.


\section{Picard modular groups}\label{arithmetic}



\subsection{Arithmetic subgroups of $\SU(2, 1)$}\label{algebraic groups}


Let $k$ be an imaginary quadratic field and $h$ a hermitian form on $k^3$ of signature $(2,1)$. If $\mc{G}_h$ is the $\Q$-algebraic group such that
\[
\mc{G}_h(\Q) \cong \{ A \in \SL_3(k)\ :\ {}^t \conj{A} h A = h \},
\]
then $\mc{G}_h(\R) \cong \SU(2,1)$, i.e., $\mc{G}_h$ is a $\Q$-form of $\SU(2,1)$. Though there are two isomorphism classes of hermitian forms over $k$ with given signature and dimension, all the $\mc{G}_h$ over a fixed $k$ are $\Q$-isomorphic \cite[$\S$1.2]{Prasad--Yeung}.

The subgroup
\[
\mc{G}_h(\Z) = \{ A \in \SL_3(\mc{O}_k)\ :\ {}^t \conj{A} h A = h \},
\]
is a nonuniform lattice in $\SU(2,1)$. Its commensurability class depends only on $k$, and these lattices give all nonuniform arithmetic lattices in $\SU(2,1)$ up to commensurability. That is, commensurability classes of nonuniform arithmetic lattices are in one-to-one correspondence with imaginary quadratic fields. Any lattice commensurable with $\mc{G}_h(\Z)$ is called a \emph{Picard modular group}. Special amongst the Picard modular groups are the \emph{principal arithmetic lattices}, which we now describe.


\subsection{Principal arithmetic lattices}\label{arithmetic lattices}


Fix an imaginary quadratic field $k$ with integer ring $\mc{O}_k$, set
\[
h_0 = \begin{mat}
0 &
0 &
1 \\

0 &
-1 &
0 \\

1 &
0 &
0
\end{mat},
\]
and let $\mc{G}$ be the associated $\Q$-form of $\SU(2,1)$. For any prime $p$ of $\Q$, consider the $\Q_p$-points $\mc{G}(\Q_p)$ of $\mc{G}$. If $p$ splits in $\mc{O}_k$, then $\mc{G}(\Q_p)$ is isomorphic to $\SL_3(\Q_p)$. When $p$ is inert or ramifies, let $k_p$ be the completion of $k$ at the prime ideal of $\mc{O}_k$ above $p$. This is a quadratic extension of $\Q_p$, and $\mc{G}(\Q_p)$ is the unitary group of $h_0$ under this extension.

For each $p$, let $K_p$ be a compact open subgroup of $\mc{G}(\Q_p)$, and suppose that $K_p = \mc{G}(\Z_p)$ for all but finitely many $p$. (Note that this integral structure comes from our above choice of basis. The important point is that $\mc{G}(\Z_p)$ is \emph{hyperspecial} for each $p$ \cite[$\S$3.5]{Tits}.) Then
\[
K_f = \prod_{p} K_p < \mc{G}(\mathbb{A}_f)
\]
is a compact open subgroup, where $\mathbb{A}_f$ is the finite adele ring of $\Q$, and $\Gam_{K_f} = K_f \cap \mc{G}(\Q)$ defines a lattice in $\SU(2,1)$ commensurable with $\mc{G}(\Z)$. In fact, $\mc{G}(\Z)$ is amongst these lattices, where $K_p = \mc{G}(\Z_p)$ for all $p$. We will always denote the lattice $\mc{G}(\Z)$, the standard Picard modular, by $\Gam_{\textrm{std}}$.

The conjugacy classes of maximal compact open subgroups of $\mc{G}(\Q_p)$ are determined by the vertices of the local Dynkin diagram of $\mc{G}(\Q_p)$. The cases of interest to us are worked out in detail in the examples of \cite{Tits}. However, a concrete set of maximal parahoric subgroups will be helpful in understanding our method of counting cusps, so we briefly recall the construction, via stabilizers of lattices in $k_p^3$, below.


\subsection{Parahoric subgroups}\label{explicit parahorics}


When $p$ splits in $\mc{O}_k$, representatives for the conjugacy classes of maximal compact open subgroups, i.e., the maximal parahoric subgroups, of $\mc{G}(\Q_p)$ are
\[
K_p^{v_0} = \SL_3(\Z_p)
\]
\[
K_p^{v_1} = \left\{ \begin{mat}
a_{1 1} &
\frac{1}{p} a_{1 2} &
\frac{1}{p} a_{1 3} \\

p a_{2 1} &
a_{2 2} &
a_{2 3} \\

p a_{3 1} &
a_{3 2} &
a_{3 3}
\end{mat} \in \SL_3(\Q_p)\ :\ a_{j k} \in \Z_p \right\}
\]
\[
K_p^{v_2} = \left\{ \begin{mat}
a_{1 1} &
a_{1 2} &
\frac{1}{p} a_{1 3} \\

a_{2 1} &
a_{2 2} &
\frac{1}{p} a_{2 3} \\

p a_{3 1} &
p a_{3 2} &
a_{3 3}
\end{mat} \in \SL_3(\Q_p)\ :\ a_{j k} \in \Z_p \right\}.
\]
If $x_1,x_2,x_3$ is the corresponding basis for $\Q_p^3$, then $K_p^{v_0}$ is the stabilizer of the standard $\Z_p$-lattice,
\[
K_p^{v_1} = \mathrm{stab} \left\{ \frac{a}{p} x_1 + b x_2 + c x_3\ :\ a,b,c \in \Z_p \right\},
\]
\[
K_p^{v_2} = \mathrm{stab} \left\{ \frac{a}{p} x_1 + \frac{b}{p} x_2 + c x_3\ :\ a,b,c \in \Z_p \right\}.
\]
Equivalently, these are the stabilizers of three vertices, $v_0$, $v_1$, and $v_2$, of the corresponding Bruhat-Tits building. These vertices determine a chamber $t_0$ of the building, which is a triangle. Let $K_p^0$ be the associated Iwahori subgroup
\[
K_p^{v_0} \cap K_p^{v_1} \cap K_p^{v_2},
\]
that is, the group of elements which are upper-triangular modulo $p$.

If $p$ is inert or ramified, the building is a tree. Let $\mc{O}_p$ be the ring of integers of the quadratic extension $k_p$ of $\Q_p$ and $\pi$ a uniformizer for $\mc{O}_p$. Extend complex conjugation on $k / \Q$ to the nontrivial Galois involution of $k_p / \Q_p$.

There are two conjugacy classes of compact open subgroups,
\[
K_p^{v_0} = \mc{G}(\Z_p) = \left\{ A \in \SL_3(\mc{O}_p)\ :\ {}^t \conj{A} h_0 A = h_0 \right\},
\]
the stabilizer in $\SU(h_0, k_p)$ of the standard $\mc{O}_p$-lattice, and
\[
K_p^{v_1} = \left\{ A = \begin{mat}
a_{1 1} &
\frac{1}{\pi} a_{1 2} &
\frac{1}{\pi} a_{1 3} \\

\pi a_{2 1} &
a_{2 2} &
a_{2 3} \\

\pi a_{3 1} &
a_{3 2} &
a_{3 3}
\end{mat} \in \SU(h_0, k_p)\ :\ a_{j k} \in \mc{O}_p \right\},
\]
which is the stabilizer in $\SU(h_0, k_p)$ of the $\mc{O}_p$-lattice
\[
\left\{ \frac{a}{\pi} x_1 + b x_2 + c x_3\ :\ a,b,c \in \mc{O}_p \right\}.
\]
The vertices $v_0$ and $v_1$ are adjoined by an edge $e_0$, which is a chamber of the building. Its stabilizer is the Iwahori subgroup
\[
K_p^0 = K_p^{v_0} \cap K_p^{v_1} = \left\{ \begin{mat}
a_{1 1} &
a_{1 2} &
a_{1 3} \\

\pi a_{2 1} &
a_{2 2} &
a_{2 3} \\

\pi a_{3 1} &
a_{3 2} &
a_{3 3}
\end{mat} \in \SU(h_0,\mc{O}_p) \right\}.
\]
In fact, one can show that preserving $h_0$ also forces the $(3,2)$-coordinate to be divisible by $\pi$, so $K_p^0$ is upper triangular modulo $\pi$.


\begin{rem}
When $2$ ramifies in $\mc{O}_k$, a slight modification is required at the prime above $2$. See \cite{Tits}.
\end{rem}



\subsection{Maximal lattices}\label{principal arithmetic lattices}


Since $\mc{G}$ is a $\Q$-form of $\SU(2,1)$, it has Strong Approximation. Therefore, for any $K_f$ as above, there exists $g \in \mc{G}(\Q)$ so that each factor of $g K_f g^{-1}$ is contained in one of the $K_p^{v_j}$ of $\S$\ref{explicit parahorics}. Furthermore, the action of $\mc{G}(\Q_p)$ on the Bruhat--Tits building at $p$ is \emph{special} \cite[$\S$3.2]{Tits}, from which it follows that if $\Gam < \mc{G}(\Q)$ is a maximal lattice, then $\Gam = \Gam_{K_f}$ for some compact open subgroup of $\mc{G}(\A_f)$, and each factor of $K_f$ stabilizes a vertex of the building. Therefore, we can, and do, assume that any $\Gam_{K_f} < \mc{G}(\Q)$ is such that each factor is contained in one of the groups $K_p^{v_j}$.

However, the $\Gam_{K_f}$ which are maximal in $\mc{G}(\Q)$ do not necessarily determine the maximal lattices in $\mc{G}(\R) \cong \SU(2,1)$. Let $\conj{\mc{G}}$ be the corresponding $\Q$-form of $\PU(2,1)$. Then every maximal arithmetic lattice in $\PU(2,1)$ defined by a hermitian form on $k^3$ is contained in $\conj{\mc{G}}(\Q)$ (see \cite[Proposition 4.2]{Platonov--Rapinchuk}), but the projection from $\mc{G}(\Q)$ to $\conj{\mc{G}}(\Q)$ is not necessarily surjective. One can describe this phenomenon explicitly via lifts to the general unitary group. The lift $g$ of some element of $\conj{\mc{G}}(\Q)$ to $\mathrm{GU}(2,1)$ will have a special unitary representative lying in $\mc{G}(\Q)$ if and only if $\det(g)$ is a cube in $k$.

However, each maximal lattice in $\SU(2,1)$ in our commensurability class is the normalizer in $\mc{G}(\conj{\Q})$ of some $\Gam_{K_f}$ \cite[Proposition 1.4]{Borel--Prasad}. We may also assume that each factor of $K_f$ is either one of the vertex stabilizers $K_p^{v_j}$, or is an Iwahori subgroup $K_p^0$. Even further, when $p$ is inert or ramified, every element of $\mc{G}(\conj{\Q}_p)$ acts by special automorphisms, so we may assume that $K_f$ is Iwahori only at split primes (see \cite[$\S$2.2]{Prasad--Yeung}).

Let $\conj{\Gam} < \mc{G}(\conj{\Q})$ be a maximal lattice, normalizing the lattice $\Gam_{K_f}$. Let $I$ be the set of (split) primes for which $K_f$ is Iwahori at $p$, so $K_f$ is one of the groups $K_p^{v_j}$ for all other $p$. Then, we have the following upper bound for $[\conj{\Gam} : \Gam_{K_f}]$, which will be of use several times throughout this paper.


\begin{prop}[See \cite{Borel--Prasad, Prasad--Yeung}]\label{normalizer index}
Let $\conj{\Gam}$ and $\Gam_{K_f}$ be as above. Then
\[
[\conj{\Gam} : \Gam_{K_f}] \leq 3^{1 + |I|} h_{k, 3},
\]
where $h_{k, 3}$ is the order of the 3-primary part of the class group of the field $k$ from which the algebraic group $\mc{G}$ is defined.
\end{prop}



\begin{pf}[Sketch]
Since the ideas behind the proof are relevant later, we give a sketch. If $k \neq \Q(\sqrt{-3})$, then the center of $\SU(2,1)$ is not contained in $\Gam$, since the center is $\zeta_3 \Id$, where $\zeta_3$ is a primitive $3^{rd}$ root of unity. When $k = \Q(\sqrt{-3})$, $\Gam_{K_f}$ is normalized in $\U(2,1)$ by an element conjugate to
\[
\begin{mat}
1 &
0 &
0 \\

0 &
\zeta_3 &
0 \\

0 &
0 &
1
\end{mat}.
\]
While this element has a unitary representative with entries in $\Q(\sqrt{-3})$, a special unitary representative has entries in $\Q(\sqrt[3]{\zeta_3})$. This accounts for the first factor of 3.

For any prime $p$ at which $K_f$ is Iwahori, there is an order 3 element of $\PGL_3(\Q_p)$ which permutes the three lattices defining a chamber of the building. In other words, a conjugate acts on the Bruhat--Tits building by an order three automorphism of the triangle with vertices $v_0, v_1, v_2$. If such an element preserves $h$, then it normalizes $\Gam_{K_f}$ but also does not have a special unitary representative with entries in $k$, hence the factor of $3^{|I|}$.

Since the image of $\conj{\Gam}$ in $\PU(2,1)$ is in $\conj{\mc{G}}(\Q)$, let $g$ be a representative for $g_0 \in \conj{\Gam}$ in the general unitary group with entries in $k$. Then
\[
{}^t \conj{g} h g = \al h
\]
for some $\al \in \Q^\times$, from which it follows that $|\det(g)|^2 = \al^3$. This lift is well-defined up to an element of $(k^\times)^3$, so we have a map
\[
\conj{\Gam} / \Gam_{K_f} \to \ker \left( \mathrm{N}_{k / \Q} : k^\times / (k^\times)^3 \to \Q^\times / (\Q^\times)^3 \right).
\]

There is a natural identification of this kernel with order three elements of the class group of $k$. Choose a non-unit element $x$ representing a nontrivial class in $k^\times / (k^\times)^3$ for which $\mathrm{N}_{k / \Q}(x) = y^3 \in \Q^\times$. Then, the principal ideal $(x)$ decomposes as a product of prime ideals $\mf{P}_1^{\al_1} \cdots \mf{P}_r^{\al_r}$, and
\[
\prod_{j = 1}^r \mf{P}_j^{\al_j} \conj{\mf{P}}_j^{\al_j} = (x) (\conj{x}) = (y)^3 = (p_1)^{3 \beta_1} \cdots (p_r)^{3 \beta_r},
\]
where $y = p_1^{\beta_1} \cdots p_r^{\beta_r}$ as a product of rational primes.

Since $p_j = \mathrm{N}_{k / \Q}(\mf{P}_j)$ for $p_j$ split or ramified, $\al_j = 3 \beta_j$ for those $j$. If $p_j$ is inert, then $(p_j) = \mf{P}_j$, so $2 \al_j = 3 \beta_j$, so $\al_j$ is still divisible by $3$. This implies that $(x)$ is a product of $3^{rd}$ powers of ideals,
\[
(x) = \mf{I}_1^3 \cdots \mf{I}_r^3 = (\mf{I}_1 \cdots \mf{I}_r)^3.
\]
Since $x$ is not a cube, $\mf{I}_x = \mf{I}_1 \cdots \mf{I}_r$ is not a principal ideal. Therefore $\mf{I}_x$ represents an element of the class group of $k$ with order dividing 3. That is, the class of $\mf{I}_x$ lies in the $3$-primary component of the ideal class group.
\end{pf}



\section{Cusps of Picard modular groups}\label{cusps}



\subsection{Cusps of standard Picard modular groups}\label{cusp set}


The first step in counting cusps of Picard modular groups is to describe the cusp set (see $\S$\ref{geometry}). The notation of $\S$\ref{arithmetic} is used throughout this section. Let $\conj{\Gam} < \mc{G}(\conj{\Q})$ be a maximal lattice which normalizes $\Gam_{K_f}$ as in $\S$\ref{principal arithmetic lattices}. Since $\SU(2,1) \to \PU(2,1)$ is central, and since the center of $\SU(2,1)$ acts trivially on $\mathbf{H}_\C^2 \cup \partial_\infty \mathbf{H}_\C^2$, there is no issue with considering lattices in $\SU(2,1)$, as opposed to their image in $\PU(2,1)$. The following lemma is known, cf.\ \cite{Chinburg--Long--Reid} and \cite[$\S$4.7]{Platonov--Rapinchuk}.


\begin{lem}\label{commensurability}
The cusp set of $\conj{\Gam}$ is independent of its commensurability class and equals $\mathbb{P}^2(k) \cap \partial_\infty \mathbf{H}_\C^2$, i.e., the set of $h_0$-isotropic lines in $k^3$.
\end{lem}


Therefore, to count the number of cusps of $\conj{\Gam}$, we will to count the $\conj{\Gam}$-equivalence classes of $h_0$-isotropic lines in $k^3$. To each $h_0$-isotropic line $\ell$ there corresponds an ideal class of $k$ as follows. For any $x \in \ell$, consider the fractional ideal
\[
I_x = \{ \al \in k\ :\ \al x \in \mc{O}_k^3 \}.
\]
The ideal class of $I_x$ depends only on $\ell$. Call this ideal class $\mathrm{cl}(\ell)$.

If $\ell$ and $\ell^\prime$ are $\gam$-equivalent for some $\gam$ in the standard Picard modular group for $h_0$, $\Gam_{\textrm{std}}$, then $\mathrm{cl}(\ell)=\mathrm{cl}(\ell^\prime)$, since $\gam$ preserves $\mc{L}_0 = \mc{O}_k^3$ and has determinant one. There is also a converse, due to Zink \cite{Zink}. We sketch the proof, highlighting the ideas used later in this paper. The key fact, which is the reason for our specification of $h_0$ and $\mc{L}_0$, is that $\mc{L}_0$ is unimodular for $h_0$.


\begin{thm}[\cite{Zink}]\label{zink1}
Let $\ell$ and $\ell^\prime$ be $h_0$-isotropic lines in $k^3$. Then $\ell$ and $\ell^\prime$ are $\Gam_{\textrm{std}}$-equivalent if and only if $\mathrm{cl}(\ell) = \mathrm{cl}(\ell^\prime)$. Given any ideal class $c$ of $k$, there exists an $h_0$-isotropic line $\ell$ with $\mathrm{cl}(\ell) = c$.
\end{thm}



\begin{pf}
Let $\phi : \ell \to \ell^\prime$ be an abstract $\mc{O}_k$-module isomorphism taking $\ell \cap \mc{L}_0$ to $\ell^\prime \cap \mc{L}_0$. Zink proves that there exists an $h_0$-preserving automorphism $\gam : \mc{L}_0 \to \mc{L}_0$ which restricts to $\phi$ as follows. 

There exist isotropic vectors $x \in \ell \cap \mc{L}_0$ and $y \in \mc{L}_0$ such that $h_0(x,y) = 1$, by unimodularity (and since $k^3$ has odd dimension). Then, we can find a decomposition of $\mc{L}_0$ as an $\mc{O}_k$-module
\[
I_\ell x \oplus I_1 y \oplus I_2 z,
\]
where $I_\ell$, $I_1$, and $I_2$ are fractional ideals and
\[
h_0(x,z) = h_0(y,z) = 0,
\]
\[
h_0(z,z) \neq 0.
\]
Note that $\mathrm{cl}(I_\ell) = \mathrm{cl}(\ell)$ and $I_1 \cong \Hom_{\mc{O}_k}(\conj{\ell}, \mc{O}_k)$, where $\conj{\ell}$ denotes $\ell$ with the complex conjugate $\mc{O}_k$-module structure.

There is a similar decomposition with respect to $\ell^\prime$ with generators $x^\prime$, $y^\prime$, and $z^\prime$. This defines an element $\gam$ of $\Gam_{\textrm{std}}$ which takes $\ell$ to $\ell^\prime$ by
\[
\gam(x) = x^\prime,\ \gam(y) = y^\prime,\ \gam(z) = z^\prime.
\]
This preserves $h_0$ by construction, so $\gam \in \Gam_{\textrm{std}}$. (Note that $\gam$ might not have determinant one, but after changing our choice of basis by a unit of $\mc{O}_k$, the determinant is one.)

To construct an isotropic line with given ideal class $c$, let $\mf{c}$ be a fractional ideal representing $c$ and $\mf{c}^\prime = \Hom_{\mc{O}_k}(\conj{\mf{c}},\mc{O}_k)$. Then the abstract rank two $\mc{O}_k$-module $\mf{c} \oplus \mf{c}^\prime$ is unimodular for the hermitian form $H$ defined by
\[
H(\mf{c},\mf{c}) = H(\mf{c}^\prime,\mf{c}^\prime) = 0,\quad H(x,x^\prime) = x^\prime(x),\quad x \in \mf{c},\ x^\prime \in \mf{c}^\prime.
\]
Set $\mf{d} = \mf{c}^{-2}$ and extend $H$ to the rank three $\mc{O}_k$-module $\mc{L} = \mf{c} \oplus \mf{c}^\prime \oplus \mf{d}$ by
\[
H(x,y) = \frac{x \conj{y}}{\mathrm{N}(\mf{c})^2},\ x,y \in \mf{d}
\]
\[
H(\mf{c},\mf{d}) = H(\mf{c}^\prime,\mf{d}) = 0.
\]

Then $\mc{L}$ is an odd unimodular lattice of rank three with respect to $H$, which has signature $(2,1)$. It has the same class as $\mc{L}_0$ by construction, so it is isomorphic to $\mc{L}_0$ equipped with $h_0$. The line generated by the image of $\mf{c}$ in $k^3$ is an isotropic line with ideal class $c$.
\end{pf}



\begin{cor}[\cite{Zink, Zeltinger}]\label{standard Picard modular group cusps}
The standard Picard modular group $\Gam_{\textrm{std}}$ has $h_k$ cusps, where $h_k$ is the class number of $k = \Q(\sqrt{-d})$. More generally, if $\mc{L}$ is a lattice in $k^3$ which is unimodular with respect to some hermitian form $h$, then the lattice $\Aut(\mc{L}, h)$ in $\SU(2,1)$ has $h_k$ cusps.
\end{cor}


The remainder of this section generalizes Zink's ideas and the general approach of \cite{Chinburg--Long--Reid} to count the number of cusps for $\conj{\Gam}$ and some of its congruence subgroups. First, we consider cusps of a suitable generalization of the Hecke congruence subgroup $\Gam_0(N)$.  This count is then used to count cusps of the maximal lattice.


\subsection{Cusps of congruence subgroups}\label{congruence subgroup cusps}


We now define the subgroups $\Gam(\mc{P}_1, \mc{P}_2, \mc{B})$ of $\Gam_{\textrm{std}}$ that generalize Hecke's subgroups $\Gam_0(N)$ of the modular group. Let $\mc{P}_1$ and $\mc{P}_2$ be disjoint finite sets of rational primes which split in $\mc{O}_k$, and $\mc{B}$ a third disjoint finite set of primes. Then $\Gam(\mc{P}_1, \mc{P}_2, \mc{B})$ is the subgroup of those $\gam \in \Gam_{\textrm{std}}$ for which the reduction of $\gam$ modulo $p$ has the form
\[
\begin{mat}
* &
* &
* \\

* &
* &
* \\

0 &
0 &
*
\end{mat} \quad \quad p \in \mc{P}_1
\]
\[
\begin{mat}
* &
* &
* \\

0 &
* &
* \\

0 &
* &
*
\end{mat} \quad \quad p \in \mc{P}_2
\]
and is upper triangular modulo $p$ for all $p \in \mc{B}$.

Note that $\mc{P}_1$ and $\mc{P}_2$ determine the two parabolic subgroups of $\SL_3(\F_p)$ containing the Borel subgroup of upper triangular matrices. The groups $\SU(3, \F_p^2)$ and $\SO(3, \F_p)$ have no proper parabolic subgroups which properly contain a Borel subgroup, which is why all the primes in $\mc{P}_1$ and $\mc{P}_2$ are split. Write $\mc{B} = \mc{B}^s \cup \mc{B}^i \cup \mc{B}^r$ as a disjoint union of its split, inert, and ramified primes, respectively.


\begin{prop}\label{picard gam_0 cusps}
The space $\Hy_\C^2 / \Gam(\mc{P}_1, \mc{P}_2, \mc{B})$ has
\[
2^{|\mc{P}_1| + |\mc{P}_2| + |\mc{B}^i| + |\mc{B}^r|} 3^{|\mc{B}^s|} h_k
\]
cusps.
\end{prop}



\begin{pf}
By Strong Approximation for $\Gam_{\textrm{std}}$, it suffices to consider one prime $p$. Since the reduction of $\Gam_{\textrm{std}}$ modulo $p$ maps onto the corresponding $\F_p$-group \cite[$\S$3.5]{Tits}, it suffices to count $\Gam(\mc{P}_1, \mc{P}_2, \mc{B})$-orbits of isotropic lines modulo $p$. Each algebraic group $\mc{G}(\F_p)$ over $\F_p$, $\SL_3(\F_p)$, $\SU(3, \F_p^2)$, or $\SO(3, \F_p)$, acts transitively on the isotropic lines modulo $p$.

For the primes in $\mc{P}_1$, $\mc{P}_2$, or a nonsplit prime of $\mc{B}$, there are two orbits of isotropic lines modulo $p$ under the corresponding reduction modulo $p$. In each case, the two orbits are represented by $(1, 0, 0)$ and $(0, 0, 1)$. For the split primes in $\mc{B}$, there are three.

For each orbit of isotropic lines modulo $p$ and each element $c$ of the class group of $k$, it remains to show that there is an isotropic line in $k^3$ with ideal class $c$ and whose reduction modulo $p$ lands the specified orbit. Let $\ell^\prime$ be any isotropic line with $\textrm{cl}(\ell^\prime) = c$ and let $\conj{\ell}^\prime$ be the reduction of $\ell \cap \mc{L}_0$ modulo $p$. Then there exists $\conj{\gam} \in \mc{G}(\F_p)$ so that $\conj{\gam}(\conj{\ell}^\prime) = \conj{\ell}$. Since the reduction is surjective, there exists $\gam \in \Gam_{\textrm{std}}$ whose reduction is $\conj{\gam}$. Then $\ell = \gam(\ell^\prime)$ has reduction modulo $p$ equal to $\conj{\ell}$. This completes the proof.
\end{pf}


The following lemma will allow us to apply Proposition \ref{picard gam_0 cusps} to study a general maximal lattice.


\begin{lem}\label{intersection is a gam_0}
Let $\conj{\Gam} < \mc{G}(\R)$ be a maximal lattice commensurable with the standard Picard modular group $\Gam_{\textrm{std}}$. Then, we can conjugate $\conj{\Gam}$ so that $\conj{\Gam} \cap \Gam_{\textrm{std}}$ is of the form $\Gam(\mc{P}_1, \mc{P}_2, \mc{B})$.
\end{lem}



\begin{pf}
Suppose that $\conj{\Gam}$ normalizes the principal arithmetic lattice $\Gam_{K_f} = \conj{\Gam} \cap \mc{G}(\Q)$. We can assume after a conjugation that $K_f$ is contained in a product of the $K_p^{v_j}$ as above, so
\[
\conj{\Gam} \cap \Gam_{\textrm{std}} = \conj{\Gam} \cap \mc{G}(\Q) \cap \Gam_{\textrm{std}} = \Gam_{K_f} \cap \Gam.
\]
The lemma follows from the description of each $K_p^{v_j}$ given in $\S$\ref{explicit parahorics}.
\end{pf}



\subsection{Proof of the main result}\label{maximal picard cusps}


Let $\conj{\Gam} < \SU(2, 1)$ be a maximal arithmetic lattice, $k$ be the number field from which it is defined, and $\Gam_{K_f} < \mc{G}(\Q)$ the principal arithmetic lattice that it normalizes, where $\mc{G}$ is the $\Q$-form of $\SU(2, 1)$ determined by $k$ and the hermitian form
\[
h_0 = \begin{mat}
0 &
0 &
1 \\

0 &
-1 &
0 \\

1 &
0 &
0
\end{mat}.
\]
Conjugate so that the $p$-component of $K_f$ equals $K_p^{v_j}$ or the Iwahori subgroup $K_p^0$ for every $p$, and is only equal to $K_p^0$ for some finite set of rational primes which split in $\mc{O}_k$.

Let $I$ be the set of primes for which the $p$-component of $K_f$ is $K_p^0$. Let $\Xi_{K_f}$ be, as in \cite{Borel--Prasad} (their notation is $\Xi_\Theta$), the set of $p$ in $I$ such that $\conj{\Gam}$ contains an element $c_p$ which acts nontrivially on the chamber of the Bruhat--Tits building at $p$ corresponding to the Iwahori subgroup $K_p^0$. Then $c_p$ normalizes $K_p^0$, and its image in $\PGL_3(\Q_p)$ has order 3.

Finally, let $\Gam_{\textrm{std}}$ be the standard Picard modular group, $\mc{G}(\Z)$ for $h_0$, and $\conj{\Gam}_{\textrm{std}}$ its normalizer in $\mc{G}(\conj{\Q})$. Now, we prove some auxiliary results necessary for the proof of Theorem \ref{intro cusp count}.


\begin{prop}\label{standard normalizer}
If $\Gam_{\textrm{std}}$ is the standard Picard modular group for $h_0$ as above, then $[\conj{\Gam}_{\textrm{std}} : \Gam_{\textrm{std}}] = 3 h_{k, 3}$.
\end{prop}



\begin{pf}
One can prove this directly via Galois cohomology and Strong Approximation, as in \cite{Margulis--Rohlfs} or \cite{Borel--Prasad}, but we give an elementary (though morally equivalent) argument. By Proposition \ref{normalizer index},
\[
[\conj{\Gam}_{\textrm{std}} : \Gam_{\textrm{std}}] \leq 3 h_{k, 3},
\]
so it suffices to prove the opposite inequality. The element of determinant $\zeta_3$ from the proof of Proposition \ref{normalizer index} generates the factor of 3 when $k = \Q(\sqrt{-3})$, and completes the proof in that case. The center of $\SU(2,1)$, thus of $\conj{\Gam}$ is cyclic of order 3, generated by the scalar matrix $\zeta_3 \Id$. When $k \neq \Q(\sqrt{-3})$, this subgroup is not in $\Gam$, so it generates a three-fold extension. This completes the proof when $h_{k, 3} = 1$. Now, we need to exhibit an element $g \in G(\conj{\Q})$ that normalizes $\Gam$ for each order three element of the class group of $k$.

Let $J \subset \mc{O}_k$ be a non-principal ideal such that $J^3$ is principal. If $\{e_1, e_2, e_3\}$ are the standard basis vectors for $k^3$, then
\[
\mc{L}_0 = \mc{O}_k^3 = \mc{O}_k e_1 \oplus \mc{O}_k e_2 \oplus \mc{O}_k e_3.
\]
Set
\[
\mc{L}_1 = J e_1 \oplus J e_2 \oplus J e_3 \subset \mc{L}_0.
\]
Then $\mc{L}_1$ is $(n)$-modular for $h$, where $(n) = J \conj{J}$. Since $J^3$ is principal, $\textrm{cl}(\mc{L}_1) = \textrm{cl}(\mc{L}_0)$, so we can find a decomposition of $\mc{L}_1$, akin to that of Theorem \ref{zink1}, as
\[
\mc{O}_k e_1^\prime \oplus \mc{O}_k e_2^\prime \oplus \mc{O}_k e_3^\prime,
\]
where $h(e_{j_1}^\prime, e_{j_2}^\prime) = 0$ or $\pm n$. This means that there is an isomorphism $g$ from $\mc{L}_0$ to $\mc{L}_1$, sending $e_j$ to $e_j^\prime$. The linear extension of $g$ to $k^3$ then, by definition, lies in the general unitary group of $h$ and scales $h$ by $n$.

We claim that $g$ normalizes $\Gam_{\textrm{std}}$. Indeed, since $J$ is an ideal of $\mc{O}_k$, it follows that any $\gam \in \Gam$ fixes $\mc{L}_1$, so $g^{-1} \gam g$ is an automorphism of $\mc{L}_0$ of determinant one. Since $g$ is in the general unitary group, $g^{-1} \gam g$ is still in the special unitary group. In particular, $g^{-1} \gam g \in \Gam_{\textrm{std}}$.

Finally, let $g_I$ and $g_J$ be two elements as above, constructed with respect to the ideals $I$ and $J$, where each represents distinct three-torsion in the class group. Then the cosets $g_I \Gam_{\textrm{std}}, g_J \Gam_{\textrm{std}} \subset \conj{\Gam}_{\textrm{std}}$ are disjoint. Indeed, $g_I$ and $g_J$ will scale $h$ by different factors, even up to multiplication by scalar matrices (which multiply determinants by a cube of $k$), so their cosets must be distinct. This proves that $[\conj{\Gam}_{\textrm{std}} : \Gam_{\textrm{std}}] \geq 3 h_{k, 3}$ and completes the proof of the proposition.
\end{pf}



\begin{cor}\label{normalizers dont fix ideal classes}
Let $\Gam_{\textrm{std}}$ be the standard Picard modular group for $h$ as above, $\conj{\Gam}_{\textrm{std}}$ its normalizer, and $g \in \conj{\Gam}_{\textrm{std}}$ one of the elements constructed in Proposition \ref{standard normalizer}, built with respect to the ideal $J$ of $\mc{O}_k$. Then $g$ acts on the cusps of $\Gam$ by
\[
\textrm{cl}(\ell) \overset{g}{\mapsto} \textrm{cl}(J)^{-1} \textrm{cl}(\ell).
\]
In particular, $\textrm{cl}(g(\ell)) \neq \textrm{cl}(\ell)$ for any isotropic line $\ell \subset k^3$.
\end{cor}



\begin{pf}
Recall that $\textrm{cl}(\ell)$ is computed with respect to $\mc{L}_0$, not $\mc{L}_1 = g(\mc{L}_0)$. Let $\ell$ be an isotropic line, $x \in \ell$, and $I_x$ the fractional ideal of $k$ associated to $\ell$, so $\textrm{cl}(\ell) = \textrm{cl}(I_x)$. Consider the ideal $I_{g(x)}$, the set of $\lam \in k$ for which $\lam g(x) \in \mc{L}_0$. If $\lam g(x) \in \mc{L}_0$, then $j \lam g(x) \in \mc{L}_1$ for all $j \in J$, so
\[
j \lam x \in g^{-1}(\mc{L}_1) = \mc{L}_0.
\]
That is, $j \lam \in I_x$ for all $j \in J$, so $J I_{g(x)} \subseteq I_x$.

Now, consider $J^{-1} I_x$. Since
\[
J^{-1} = \Hom_{\mc{O}_k}(J, \mc{O}_k),
\]
if $j^{-1} \mu \in J^{-1} I_x$,
\[
j^{-1} \mu g(x) = j^{-1} g(\mu x),
\]
and $\mu x \in \mc{L}_0$ by definition, so $g( \mu x) \in \mc{L}_1$. Write
\[
g(\mu x) = j_1 e_1 + j_2 e_2 + j_3 e_3,
\]
with $j_\al \in J$, $\al = 1,2,3$. Then $j^{-1}(j_\al) \in \mc{O}_k$ for all $\al = 1,2,3$, so $j^{-1} \mu \in I_{g(x)}$. Thus, $J^{-1} I_x \subseteq I_{g(x)}$, i.e., $I_x \subseteq J I_{g(x)}$, so the two ideals are equal. This implies that $\textrm{cl}(g(\ell)) = \textrm{cl}(J)^{-1} \textrm{cl}(\ell)$, as desired.
\end{pf}


Now we are ready to state and prove the main technical result of this paper.


\begin{thm}\label{maximal count}
Let $\conj{\Gam} < \SU(2, 1)$ be a maximal arithmetic lattice, and $\Gam_{K_f}$ the principal arithmetic lattice that it normalizes. Let $I$ be the set of (necessarily split) rational primes for which the $p$-component of $K_f$ is an Iwahori subgroup, and $m_{K_f} = |I \ssm \Xi_{K_f}|$. Then $\Hy_\C^2 / \conj{\Gam}$ has
\[
3^{m_{K_f}} \frac{h_k}{h_{k, 3}}
\]
cusps.
\end{thm}



\begin{pf}
Let $\Gam_{\textrm{std}}$ be the standard Picard modular group for $h_0$ as above (so Theorem \ref{zink1} applies) and $\conj{\Gam}_{\textrm{std}}$ its normalizer. Consider the diagram of subgroups in Figure \ref{diagram}, from which we will prove the theorem two stages.

\begin{figure}[htb]
\begin{center}
\[
\xymatrix{
\conj{\Gam} \ar@{-}[dr] \ar@{-}[d] &
 &
\conj{\Gam}_{\textrm{std}} \ar@{-}[d] \ar@{-}[dl] \\
\Gam_{K_f} \ar@{-}[dr] &
\conj{\Gam} \cap \conj{\Gam}_{\textrm{std}} \ar@{-}[d] &
\Gam_{\textrm{std}} \ar@{-}[dl] \\
 &
\Gam(\mc{P}_1, \mc{P}_2, \mc{B}) &
 \\
}
\]
\caption{The intersection of $\conj{\Gam}$ and $\conj{\Gam}_{\textrm{std}}$.}
\label{diagram}
\end{center}
\end{figure}

First, assume $k \neq \Q(\sqrt{-3})$. We will deal with $\Q(\sqrt{-3})$ at the end of the proof. Then $\conj{\Gam} \cap \Gam_{\textrm{std}} = \Gam(\mc{P}_1, \mc{P}_2, \mc{B})$ has $2^{m_1} 3^{m_2} h_k$ cusps by Lemma \ref{picard gam_0 cusps}, where $m_1$ and $m_2$ are explicitly determined by $K_f$. We claim that if $\ell$ is an isotropic line and $g \in \conj{\Gam} \cap \conj{\Gam}_{\textrm{std}}$ is such that $\textrm{cl}(\ell) = \textrm{cl}(g(\ell))$, then either $g \in \Gam(\mc{P}_1, \mc{P}_2, \mc{B})$ or $g$ is in the center of $\SU(2,1)$, which is cyclic of order three and acts trivially on all cusps. Equivalently, $\conj{\Gam} \cap \conj{\Gam}_{\textrm{std}}$ has
\[
\frac{2^{m_1} 3^{m_2} h_k}{\frac{1}{3}[\conj{\Gam} \cap \conj{\Gam}_{\textrm{std}} : \Gam(\mc{P}_1, \mc{P}_2, \mc{B})]}
\]
cusps.

If $g$ is in the center, $\zeta_3 \notin \Q(\sqrt{-d})$ since $d \neq 3$, so $g \notin \Gam$. That the center acts trivially on isotropic lines is evident. Now, suppose $g$ is not central and that $\textrm{cl}(\ell) = \textrm{cl}(g(\ell))$ for $g \in \conj{\Gam} \cap \conj{\Gam}_{\textrm{std}}$. We must show that $g \in \Gam_{\textrm{std}}$, which is precisely the content of Corollary \ref{normalizers dont fix ideal classes}.

Now, consider $\conj{\Gam} \cap \conj{\Gam}_{\textrm{std}} < \conj{\Gam}$. First, this contains elements $c_p = c$ for $p \in I$, where $I$ is the set of (split) primes for which the $p$-component of $K_f$ is $K_p^0$. These elements generate the extension $\conj{\Gam} / \Gam_{K_f}$, and $\conj{\Gam} / \Gam_{K_f}$ is an elementary abelian 3-group. Every such element identifies the three equivalence classes of cusps (for each ideal class) corresponding to the fact that the reduction modulo $p$ of $\Gam_{K_f} \cap \Gam_{\textrm{std}}$ is the Borel subgroup of $\SL_3(\F_p)$. There are, by definition, $\Xi_{K_f}$ such elements.

Also, for each $p \notin I$ such that $K_p$ is not $G(\Z_p) = K_p^{v_0}$, $\Gam_{K_f}$ contains an element which identifies the two classes of isotropic lines modulo $p$ which are not identified by $\conj{\Gam} \cap \conj{\Gam}_{\textrm{std}}$. One can see this explicitly via Strong Approximation and the reduction modulo $p$ of $\Gam_{K_f}$ as described in \cite[$\S$3.11]{Tits}. This shows that $\conj{\Gam}$ has
\[
\frac{3^{|I \ssm \Xi_{K_f}|} h_k}{\frac{1}{3}[\conj{\Gam} \cap \conj{\Gam}_{\textrm{std}} : \Gam(\mc{P}_1, \mc{P}_2, \mc{B})]}
\]
cusps. It remains to show that this equals the expression in the statement of the theorem.

To show that $[\conj{\Gam} \cap \conj{\Gam}_{\textrm{std}} : \Gam(\mc{P}_1, \mc{P}_2, \mc{B})] = 3 h_{k, 3}$, it suffices by Proposition \ref{standard normalizer} to show that
\[
[\conj{\Gam}_{\textrm{std}} : \conj{\Gam} \cap \conj{\Gam}_{\textrm{std}}] = [\Gam_{\textrm{std}} : \Gam(\mc{P}_1, \mc{P}_2, \mc{B})].
\]
These intersections are exactly determined by reductions modulo $p$, where the reduction of $\conj{\Gam}_{\textrm{std}}$ modulo $p$ is a subgroup of the general linear, unitary, or orthogonal group modulo $p$, depending on the decomposition of $p$ in $\mc{O}_k$. This is well-defined by realizing of each element of $\conj{\Gam}_{\textrm{std}}$ as an element of the general unitary group of $h$ with entries in $\mc{O}_k$.

Then, the reduction of $\conj{\Gam} \cap \conj{\Gam}_{\textrm{std}}$ modulo $p$ is the appropriate parabolic or Borel subgroup of the general linear, unitary, or orthogonal group. In each case, the Borel subgroup, and thus each parabolic subgroup, contains the center. This implies that the index of the parabolic subgroup in the full $\F_p$ group is equal to the corresponding index for the reduction modulo $p$ of $\Gam_{\textrm{std}}$. This means precisely that $[\conj{\Gam}_{\textrm{std}} : \conj{\Gam} \cap \conj{\Gam}_{\textrm{std}}]$ equals $[\Gam_{\textrm{std}} : \Gam(\mc{P}_1, \mc{P}_2, \mc{B})]$, as required. This proves the theorem for $d \neq 3$.

Now, consider the case $k = \Q(\sqrt{-3})$. The group $\conj{\Gam} \cap \conj{\Gam}_{\textrm{std}}$ is generated over $\Gam(\mc{P}_1, \mc{P}_2, \mc{B})$ by a special unitary representative of the element of determinant $\zeta_3$ from the proof of Proposition \ref{standard normalizer}, which acts trivially on cusps. This proves that $\conj{\Gam} \cap \conj{\Gam}_{\textrm{std}}$ has
\[
\frac{2^{m_1} 3^{m_2} h_k}{\frac{1}{3}[\conj{\Gam} \cap \conj{\Gam}_{\textrm{std}} : \Gam(\mc{P}_1, \mc{P}_2, \mc{B})]}
\]
cusps, as in the other cases. The proof going from $\conj{\Gam} \cap \conj{\Gam}_{\textrm{std}}$ to $\conj{\Gam}$ is exactly the same as above.
\end{pf}



\subsection{Commensurability classes containing $N$-cusped elements}\label{cusp corollaries}



\begin{thm}\label{one cusp orbifolds}

For any natural number $N$, there only exist finitely many commensurability classes of arithmetic complex hyperbolic 2-orbifolds containing an orbifold with at most $N$ cusps.

\end{thm}



\begin{pf}

It suffices to prove that
\[
\frac{h_k}{h_{k,3}} \to \infty
\]
as $d \to \infty$, where $k=\Q(\sqrt{-d})$. It is a classical result of Siegel that for every $\ep > 0$ there is a constant $c(\ep)$ such that
\[
h_k > c(\ep) |d_k|^{\frac{1}{2} - \ep},
\]
where $d_k$ is the discriminant. More recently, Ellenberg and Venkatesh \cite{Ellenberg--Venkatesh} proved that there is another positive constant $c^\prime(\ep)$ so that
\[
h_{k,3} \leq c^\prime(\ep) |d_k|^{\frac{1}{3} + \ep}.
\]
(Also, see \cite{Pierce} and \cite{Helfgott--Venkatesh} for earlier bounds which suffice here.) Therefore
\[
\frac{h_k}{h_{k, 3}} \geq \frac{c(\ep)}{c^\prime(\ep)} |d_k|^{\frac{1}{6} - 2 \ep}.
\]
Fixing a small $\ep$ and letting $d \to \infty$ proves the theorem, since
\[
d_k = \left\{ \begin{matrix}
d &
d \equiv 1 \pmod{4} \\

4d &
d \not\equiv 1 \pmod{4}.
\end{matrix}\right.
\]
for square-free $d$.

\end{pf}



\begin{rem}

It is relatively easy to show that there are infinitely many distinct one-cusped elements in each commensurability class for which $h_{k, 3} = h_k$. Choosing $K_f$s for which the non-hyperspecial vertex is chosen at arbitrarily large inert primes gives an infinite family of non-isomorphic lattices, all of which are 1-cusped. This construction works for any natural number $N$ for which $N = h_k / h_{k, 3}$ for some imaginary quadratic field $k$.

\end{rem}


Finally, we end with a table of the known $d$ for which $\Q(\sqrt{-d})$ determines a 1-cusped orbifold. Though the list of imaginary quadratic fields for which $h_{k, 3} = h_k$ is finite, we do not know if this table is complete (see \cite{Heath-Brown}).


\begin{table}[h]
\begin{center}
\vskip 10pt
\begin{tabular}{|c|c|}
\hline
$h_k$ & $d$ \\
\hline \hline
  & 3, 4, 7 \\
1 & 8, 11, 19 \\
  & 43, 67, 163 \\
\hline
  & 23, 31, 59, 83, 107, 139 \\
3 & 211, 283, 307, 331, 379 \\
  & 499, 547, 643, 883, 907 \\
\hline
9 & 4027 \\
\hline
27 & $\emptyset$ \\
\hline
\end{tabular} \\[10pt]
\caption{Some $\Q(\sqrt{-d})$ which give 1-cusped orbifolds.}
\vskip -20pt
\end{center}
\end{table}



\section{Higher rank}\label{su(r+1,r) lattices}



\subsection{Cusps in higher rank}\label{higher rank cusps}


We now consider nonuniform lattices in $\SU(r + 1, r)$ for $r > 1$. By Margulis's Arithmeticity Theorem, all such lattices are arithmetic, so each nonuniform arithmetic lattice in $\SU(r + 1, r)$ is commensurable with $\mc{G}(\Z)$, where $\mc{G}$ is an isotropic $\Q$-form of $\SU(r + 1, r)$. We now briefly describe the possible nonuniform commensurability classes, and describe how one determines the cusps for lattices of simple type.

In short, commensurability classes of nonuniform lattices arise from hermitian forms on $D^m$, where $D$ is a central simple division algebra over the imaginary quadratic field $k$ with involution of second kind. If $D$ has degree $d$, then $d m = 2 r + 1$. If $q = 2 r + 1$ is prime, there is a unique commensurability class of nonuniform lattices in $\SU(r + 1, r)$ for every imaginary quadratic field, and these are all the nonuniform commensurability classes. The case $d = 2 r + 1$ and $m = 1$ determines cocompact lattices, since $D$ is a division algebra.

The situation for the commensurability class of a nonuniform lattice of simple type is exactly the same as for $\SU(2,1)$. The maximal lattice $\conj{\Gam}$ normalizes a principal arithmetic lattice $\Gam_{K_f}$ for $K_f < \mc{G}(\mathbb{A}_f)$ as before. Furthermore, we may again assume that $K_f$ is hyperspecial at all but finitely many places, stabilizes a vertex of the Bruhat--Tits building at all but split primes of $\mc{O}_k$, and is an Iwahori subgroup at some finite set $I$ of split primes. In fact, the proofs carry through verbatim, with the odd integer $q$ instead of 3. See \cite{Prasad--Yeung2} for details.

In this setting, a cusp of $\Gam$ corresponds to conjugacy classes of chains of parabolic subgroups in $\Gam$. See \cite{Borel--Ji}. Recall that an \emph{isotropic flag}\index{Isotropic flag} for the hermitian form $h$ is a strictly ascending chain of $h$-isotropic subspaces
\[
\{ 0 \} = \mc{F}_0 \subset \mc{F}_1 \subset \mc{F}_2 \subset \cdots \subset \mc{F}_n \subseteq \C^n.
\]
If $h$ is a hermitian form of signature $(r + 1, r)$ on $k^q$, $q = 2 r + 1$, then a maximal isotropic flag in $k^q$ is of the form $\{ \mc{F}_j \}$, where $\mc{F}_j$ has dimension $j$, and $0 \leq j \leq r$. That is, the maximal isotropic subspace of $k^q$ has dimension $r$. For our purposes, it suffices to record the following.


\begin{prop}\label{su(r+1,r) cusps}
Suppose there is a $\Q$-form of $\SU(r + 1, r)$ of simple type over the imaginary quadratic field $k$ such that $\Gam$ is commensurable with $\mathcal{G}(\Z)$. Then cusps of $\Gam$ are in one-to-one correspondence of $\Gam$-orbits of maximal $h$-isotropic flags
\[
\mc{F} = \{ \mc{F}_j \}_{j = 0}^r \subset k^q.
\]
\end{prop}



\subsection{Zink's Theorem in higher rank}\label{higher rank zink}


We now describe Zink's Theorem for $\SU(r + 1, r)$, which is the first step to generalizing Theorem \ref{intro cusp count} to Theorem \ref{intro higher rank cusp count}. Choose $h$ to be the unimodular (for the standard $\mc{O}_k$-lattice $\mc{L}_0$) hermitian form
\[
\begin{mat}
 &
 &
 &
 &
 &
 &
1 \\

 &
 &
 &
 &
 &
\iddots &
 \\

 &
 &
 &
 &
1 &
 &
 \\

 &
 &
 &
\pm 1 &
 &
 &
 \\

 &
 &
1 &
 &
 &
 &
 \\

 &
\iddots &
 &
 &
 &
 &
 \\

1 &
 &
 &
 &
 &
 &

\end{mat},
\]
where $\pm 1$ is chosen such that $\det(h) = 1$, let $\mc{G}$ be the associated $\Q$-algebraic group, and $\Gam_{\textrm{std}} = \mc{G}(\Z)$ the special automorphism group of the hermitian lattice $(\mc{L}_0, h)$.


\begin{thm}[\cite{Zink}]\label{zink2}
Let $\Gam_{\textrm{std}}$ be as above. Then $\Gam_{\textrm{std}}$ has $h_k^r$ cusps.
\end{thm}



\begin{pf}
The proof is a direct generalization of the methods employed in the proof of Theorem \ref{zink1}. The goal is to count $\Gam_{\textrm{std}}$-orbits of maximal isotropic flags $\mc{F}$ in $k^q$. The maximal isotropic flags have dimension $r$, so a flag is a chain of subspaces
\[
F_0 = \{0\} \subset F_1 \subset \cdots \subset F_r,
\]
where $F_j$ has dimension $j$ over $k$. Therefore, $F_r \cap \mc{L}_0$ is an $\mc{O}_k$-module of rank $r$, so $F_{j + 1} / F_j$ has rank one and we can associate to $F_{j + 1}$ the ideal class of
\[
(F_{j + 1} \cap \mc{L}_0) / (F_j \cap \mc{L}_0).
\]
This associates to each flag an element of $h_k^r$.

Since $n$ is odd and $h$ is unimodular, we can find $h$-isotropic vectors $x_1, \dots x_r$ and $y_1, \dots, y_r$ so that
\[
\mc{L}_0 = I_1 x_1 \oplus \cdots \oplus I_r x_r \oplus I^\prime_1 y_1 \oplus \cdots \oplus I^\prime_r y_r \oplus I_0 z_0,
\]
where $I_j = F_j \cap \mc{L}_0$, and the $I^\prime_j$ correspond to another isotropic flag $\mc{F}^\prime$ which determines the $\mc{O}_k$-homomorphisms of $\mc{F}$ to $\mc{O}_k$. The vector $z_0$ is $h$-independent from $\mc{F}$ and $\mc{F}^\prime$ and non-isotropic. Note that this decomposition is analogous to the one used in Theorem \ref{zink1}, with `isotropic line' replaced with `isotropic flag'.

The same exact methods as in the proof of Theorem \ref{zink1} show that two isotropic flags are $\Gam$-equivalent if and only if they represent the same element of $h_k^r$, and that for each $r$-tuple in $h_k^r$ there exists a maximal isotropic flag representing that element.
\end{pf}



\subsection{Cusps of maximal lattices}\label{general proof}


We now describe how Theorem \ref{maximal count} generalizes to give Theorem \ref{intro higher rank cusp count}. Since the notions of congruence subgroups and the `Hecke congruence subgroups' $\Gam_0(N)$ become extremely complicated, given that the number of choices of proper parabolic subgroups grows with $r$, it is more natural to go straight to the main result. Recall the notation of $\S$\ref{higher rank zink}, and we assume $r > 1$, so arithmeticity is automatic.


\begin{thm}\label{general main theorem}
Let $\conj{\Gam} < \SU(r + 1, r)$ be a maximal lattice of simple type, defined via the imaginary quadratic field $k$. Conjugate $\conj{\Gam}$ so that it normalizes the principal arithmetic lattice $\Gam_{K_f}$. Let $I$ be the number of (split) primes for which $K_f$ is an Iwahori subgroup and $m_{K_f} = |I \ssm \Xi_{K_f}|$. Then $\Gam$ has
\[
q^{m_{K_f}} \frac{h_k^r}{h_{k, q}}
\]
cusps.
\end{thm}



\begin{pf}
Consider Figure \ref{diagram}, where the bottom group is now best called just $\conj{\Gam} \cap \Gam_{\textrm{std}}$. Since the $p$-component of $K_f$ is either $K_p^0$ or one of the groups $K_p^{v_j}$, each reduction modulo $p$ of $\conj{\Gam} \cap \Gam_{\textrm{std}}$ is a parabolic subgroup of the corresponding $\F_p$-group that contains the Borel subgroup of upper triangular matrices. Each parabolic subgroup has $m_p$ orbits of maximal isotropic flags modulo $p$, where $m_p = 1$ for almost every $p$. Therefore, $\conj{\Gam} \cap \Gam_{\textrm{std}}$ has
\[
h_k^r \prod_p m_p
\]
cusps. Indeed, the number of orbits of lines modulo $p$ is explicitly determined by the type of the corresponding parabolic subgroup.

Since $q$ is odd, every unit of $\mc{O}_k$ is a $q^{th}$ power, unless $q \equiv 3 \pmod{4}$ and $k = \Q(\sqrt{-3})$. We assume for now that every unit is a $q^{th}$ power in $k$. The necessary adjustments for $d = 3$ and $q \equiv 3 \pmod{4}$ are exactly the same as the case $q = 3$. Let $\conj{\Gam}_{\textrm{std}}$ be the normalizer of $\Gam_{\textrm{std}}$ in $\SU(r + 1, r)$. The center of $\SU(r + 1, r)$ is cyclic of order $q$, and $k$ contains no $q^{th}$ roots of unity for $r > 1$. The argument in Proposition \ref{standard normalizer} carries over unaltered to $q$-dimensions to show that
\[
[\conj{\Gam} : \Gam] = q h_{k, q}.
\]
See also \cite{Prasad--Yeung2}.

Furthermore, a direct generalization of Corollary \ref{normalizers dont fix ideal classes} shows that $\conj{\Gam} \cap \conj{\Gam}_{\textrm{std}}$ has
\[
\frac{h_k^r}{h_{k, q}} \prod_p m_p
\]
cusps. Indeed,
\[
[\conj{\Gam} \cap \conj{\Gam}_{\textrm{std}} : \conj{\Gam} \cap \Gam_{\textrm{std}}] = [\conj{\Gam}_{\textrm{std}} : \Gam_{\textrm{std}}],
\]
again from considering reductions modulo $p$, and each element of $\conj{\Gam}_{\textrm{std}}$ not in $\Gam_{\textrm{std}}$ acts on $h_k^r$ by multiplication by $\textrm{cl}(J)$, where $J$ is some fractional ideal for which $J^q$ is principal. (Note that $J^s$ may be principal for some divisor $s$ of $q$.)

However, for each $p$ where the $p$-component is a vertex stabilizer of the building, i.e., is not Iwahori, there exists an element of $\Gam_{K_f} < \conj{\Gam}$ for which the reduction modulo $p$ identifies the $m_p$ isotropic subspaces left inequivalent by the reduction of $\conj{\Gam} \cap \Gam_{\textrm{std}}$ modulo $p$. Since the case where $m_p = q$ corresponds precisely to when the reduction modulo $p$ of $\conj{\Gam} \cap \Gam_{\textrm{std}}$ is upper triangular, we have that $\conj{\Gam}$ has
\[
q^{m_{K_f}} \frac{h_k^r}{h_{k, q}}
\]
cusps.
\end{pf}


Suppose that $q$ is prime. Then every nonuniform arithmetic lattice in $\SU(r + 1, r)$ for $r > 1$ is of simple type. Theorem \ref{intro higher rank commensurability} is trivial for $N = 1$. Assuming the Generalized Riemann Hypothesis, Ellenberg--Venkatesh \cite{Ellenberg--Venkatesh} prove that $h_{k, q}$ grows like $|d|^{\frac{1}{2} - \frac{1}{2 q} + \ep}$. We then get Theorem \ref{intro higher rank commensurability} exactly the same as Theorem \ref{one cusp orbifolds}. Constructing infinitely many commensurable one-cusped manifolds is also done exactly the same as in the remark of $\S$\ref{cusp corollaries}. However, now the commensurability classes that produce a one-cusped element come precisely from the nine imaginary quadratic fields of class number one.


\subsection*{Acknowledgments}

This paper is a condensed version of my Ph.D.\ thesis at the University of Texas at Austin. I want to thank my advisor, Alan Reid, for his encouragement and inspiration. I also thank Jordan Ellenberg for a very helpful correspondence regarding torsion in class groups. Finally, I want to thank the people that showed interest in this project, especially Gopal Prasad, Benson Farb, Daniel Allcock, Ben McReynolds, and Misha Belolipetsky.



\bibliography{picardcusps3}

\begin{thebibliography}{10}

\bibitem{Borel--Ji}
Armand Borel and Lizhen Ji.
\newblock {\em Compactifications of symmetric and locally symmetric spaces}.
\newblock Mathematics: Theory \& Applications. Birkh\"auser Boston Inc.,
  Boston, MA, 2006.

\bibitem{Borel--Prasad}
Armand Borel and Gopal Prasad.
\newblock Finiteness theorems for discrete subgroups of bounded covolume in
  semi-simple groups.
\newblock {\em Inst. Hautes \'Etudes Sci. Publ. Math.}, (69):119--171, 1989.

\bibitem{Chinburg--Long--Reid}
Ted Chinburg, Darren Long, and Alan~W. Reid.
\newblock Cusps of minimal non-compact arithmetic hyperbolic 3-orbifolds.
\newblock {\em Pure Appl. Math. Q.}, 4(4, part 1):1013--1031, 2008.

\bibitem{Ellenberg--Venkatesh}
Jordan~S. Ellenberg and Akshay Venkatesh.
\newblock Reflection principles and bounds for class group torsion.
\newblock {\em Int. Math. Res. Not. IMRN}, (1):Art. ID rnm002, 18, 2007.

\bibitem{Heath-Brown}
D.~R. Heath-Brown.
\newblock Imaginary quadratic fields with class group exponent 5.
\newblock {\em Forum Math.}, 20(2):275--283, 2008.

\bibitem{Helfgott--Venkatesh}
H.~A. Helfgott and A.~Venkatesh.
\newblock Integral points on elliptic curves and 3-torsion in class groups.
\newblock {\em J. Amer. Math. Soc.}, 19(3):527--550, 2006.

\bibitem{Kamishima}
Yoshinobu Kamishima.
\newblock Nonexistence of cusp cross-section of one-cusped complete complex
  hyperbolic manifolds. {II}.
\newblock {\em Int. Math. Forum}, 2(25-28):1251--1258, 2007.

\bibitem{Margulis--Rohlfs}
G.~A. Margulis and J.~Rohlfs.
\newblock On the proportionality of covolumes of discrete subgroups.
\newblock {\em Math. Ann.}, 275(2):197--205, 1986.

\bibitem{Petersen}
Kathleen~L. Petersen.
\newblock One-cusped congruence subgroups of {B}ianchi groups.
\newblock {\em Math. Ann.}, 338(2):249--282, 2007.

\bibitem{Pierce}
L.~B. Pierce.
\newblock The 3-part of class numbers of quadratic fields.
\newblock {\em J. London Math. Soc. (2)}, 71(3):579--598, 2005.

\bibitem{Platonov--Rapinchuk}
Vladimir Platonov and Andrei Rapinchuk.
\newblock {\em Algebraic groups and number theory}, volume 139 of {\em Pure and
  Applied Mathematics}.
\newblock Academic Press Inc., 1994.

\bibitem{Prasad--Yeung}
Gopal Prasad and Sai-Kee Yeung.
\newblock Fake projective planes.
\newblock {\em Invent. Math.}, 168(2):321--370, 2007.

\bibitem{Prasad--Yeung2}
Gopal Prasad and Sai-Kee Yeung.
\newblock Arithmetic fake projective spaces and arithmetic fake
  {G}rassmannians.
\newblock {\em Amer. J. Math.}, 131(2):379--407, 2009.

\bibitem{Reid}
Alan~W. Reid.
\newblock Arithmeticity of knot complements.
\newblock {\em J. London Math. Soc. (2)}, 43(1):171--184, 1991.

\bibitem{Scharlau}
Winfried Scharlau.
\newblock {\em Quadratic and {H}ermitian forms}, volume 270 of {\em Grundlehren
  der Mathematischen Wissenschaften}.
\newblock Springer-Verlag, 1985.

\bibitem{Tits}
J.~Tits.
\newblock Reductive groups over local fields.
\newblock In {\em Automorphic forms, representations and {$L$}-functions
  ({P}roc. {S}ympos. {P}ure {M}ath., {O}regon {S}tate {U}niv., {C}orvallis,
  {O}re., 1977), {P}art 1}, Proc. Sympos. Pure Math., XXXIII, pages 29--69.
  Amer. Math. Soc., 1979.

\bibitem{Zeltinger}
Hanns Zeltinger.
\newblock {\em Spitzenanzahlen und {V}olumina {P}icardscher
  {M}odulvariet\"aten}.
\newblock Bonner Mathematische Schriften, 136. Universit\"at Bonn
  Mathematisches Institut, 1981.

\bibitem{Zink}
Thomas Zink.
\newblock \"{U}ber die {A}nzahl der {S}pitzen einiger arithmetischer
  {U}ntergruppen unit\"arer {G}ruppen.
\newblock {\em Math. Nachr.}, 89:315--320, 1979.

\end{thebibliography}


\end{document}